\documentclass{amsart}[12pt] 
\title[stable rank and continuous fields]{On the stable rank of algebras of 
operator fields over an n-cube}
\author{Ping Wong Ng}
\address{Department of Mathematics \\ University of Toronto \\
100 St. George St., Room 4072 \\
Toronto, Ontario \\
M5S 3G3\\
Canada}
\email{png@math.toronto.edu}

\author{Takahiro Sudo}
\address{Department of Mathematical Sciences \\ Faculty of Science \\
University of the Ryukyus \\  
Nishihara-cho, Okinawa \\ 
903-0213\\
Japan}
\email{sudo@math.u-ryukyu.ac.jp}

\date{\today}

\newtheorem{thm}{Theorem}[section]
 
\newtheorem{lem}[thm]{Lemma}
\newtheorem{cor}[thm]{Corollary}

\theoremstyle{remark}

\newcommand{\B}{\mathcal{B}}
\newcommand{\F}{\mathcal{F}}
\newcommand{\A}{\mathcal{A}}

\newcommand{\T}{\mathbb T}
\newcommand{\C}{\mathbb C}
\newcommand{\M}{\mathbb M}
\newcommand{\Z}{\mathbb Z}

\begin{document}

\begin{abstract}
Let $\A$ be a unital maximal full algebra of operator fields with 
base space $[0,1]^k$ and fibre algebras $\{ \A_t \}_{t \in [0,1]^k}$.  We 
show that the stable rank of $\A$ is bounded above by the quantity
$sup_{t \in [0,1]^k} sr( C([0,1]^k) \otimes \A_t)$.  Here the symbol
``sr" means stable rank.  Using the above estimate, we compute the stable
ranks of the $C^*$-algebras of the
(possibly higher rank) discrete Heisenberg groups.
\end{abstract}

\maketitle

\section{Introduction}

   Rieffel [\textbf{12}] introduced the notion of 
\emph{stable rank} for $C^*$-algebras
as the noncommutative
version of complex dimension of ordinary topological spaces.  It turns out 
that
the stable rank of a unital $C^*$-algebra is the same as its Bass stable
rank (see [\textbf{7}]).  The purpose of this 
paper is to study stable rank for continuous field $C^*$-algebras.
Our main result is
\begin{thm}
Let $\A$ be a unital maximal full algebra of operator fields with base space
$[0,1]^k$ and fibre algebras $\{ \A_t \}_{ t \in [0,1]^k }$.  Then the 
stable
rank of $\A$ satisfies the following inequality:
\[
sr( \A ) \leq sup_{ t \in [0,1]^k } sr ( C( [0,1]^k ) \otimes \A_t ).
\]
\end{thm}
\noindent  One of the key ingredients of the proof of the above result is
Nistor's notion of \emph{absolute connected stable rank} (see [\textbf{10}]).  
Recall
that for a unital $C^*$-algebra $\A$, the absoluted connected stable rank of
$\A$ is numerically the same as the stable rank of the tensor product
$C[0,1] \otimes \A$.  

     As an application of our theorem, we will compute the stable ranks of the
universal $C^*$-algebras of the (possibly higher rank)
discrete Heisenberg groups.  Recall that for a positive integer $n$, the
discrete Heisenberg group of rank $2n + 1$ is the group of all $2n+1$ by
$2n + 1$ upper triangular matrices with integer entries,
with ones on the diagonal, and zero 
entries on all but the first row, last column and the diagonal.  The 
rank $2n + 1$ discrete Heisenberg group is naturally a lattice subgroup of
the $2n+1$-dimensional Heisenberg Lie group.   
  
     The stable ranks of the universal $C^*$-algebras of various type I 
Lie groups
have been extensively studied (see [\textbf{13}], [\textbf{15}], [\textbf{16}]).  In particular, we point out that  
for a simply connected nilpotent Lie group $G$, the stable rank of the     
universal $C^*$-algebra $C^*(G)$ of $G$ has been computed by Sudo and Takai
(see [\textbf{15}]).  Roughly speaking, they showed that the stable rank
of $C^*(G)$ is controlled by the ordinary topological dimension of the 
space of one-dimensional representations of $G$. 

     Recently, the stable ranks of a class of non-type I solvable 
Lie groups (which
include the Mautner group) have also been computed (see [\textbf{14}]). 

     Our computations for the discrete Heisenberg groups constitute
a class of interesting nontrivial examples of the stable rank of the 
universal $C^*$-algebra
of a non-type I discrete group. 

     In a later paper,  we will apply the techniques 
developped in this paper to compute the stable ranks of arbitrary 
finitely generated, torsion-free two-step nilpotent groups. 

     A  general reference for stable rank of $C^*$-algebras is [\textbf{12}].   
General references for algebras of operator fields are [\textbf{5}],  
[\textbf{8}] and
[\textbf{17}].
 
     In what follows, for a $C^*$-algebra $\A$, the notation ``$sr(\A)$" will
always mean the stable rank of $\A$.
If, in addition, $\A$ is unital and $N$ is a positive integer, then  
$Lg_N(\A)$ is the set of $N$-tuples $(a_1, a_2, ..., a_N)$ in
$\A^N$ such that $\sum_{i = 1}^N (a_i)^* a_i$ is an invertible element in 
$\A$. 

\section{Main results}

\begin{lem} Let $\A$ be a unital maximal
full algebra of operator fields with base 
space $[0,1]$ and fibre algebras $\{ \A_t \}_{t \in \mathcal{F}}$. Then
the stable rank of $\A$ satisfies the inequality
\[
sr(\A) \leq sup_{ t \in [0,1] } sr( C[0,1] \otimes \A_t ).
\] 
\end{lem}
\begin{proof}

     In this proof, we will always let ``$\mathcal{F}$" denote the continuity 
structure for the continuous field decomposition of $\A$ in the hypothesis. 

   Suppose that $M = sup_{t \in [0,1]} sr ( C[0,1] \otimes \A_t )$ is a finite
number.  Let an 
$M$-tuple $( a_1, a_2, ... , a_M ) \in \A^M$ and a positive real number 
$\epsilon > 0$ be given.  Since $[0,1]$ is compact, let
 $I_1, I_2, ... ,I_N$ be a finite set of open intervals which cover 
$[0,1]$ and for each $i = 1, 2 ,... ,N$, let $( f_{i,1}, f_{i,2}, ... ,       
f_{i,M} )$ be an $M$-tuple in ${\A}^M$ such that 
\begin{enumerate}
\item  $\sum_{j=1}^M f_{i,j}(t)^* f_{i,j}(t)$ is an invertible element 
of $\A_t$
for $t \in I_i$.  Here the operator field $t \mapsto f_{i,j}(t)$ for $t \in 
[0,1]$, is the representation of $f_{i,j}$ as a continuous operator field
in $\A$.  And  
\item $ \| f_{i,j} - a_j \| < \epsilon$, for all $i$, $j$. 
\end{enumerate}
     We can choose such intervals $I_i$ and such elements $f_{i,j}$, since 
our hypothesis implies that the stable rank of each fibre $\A_t$ is less 
than or equal to $M$, and by the existence and continuity of operator fields 
in a  full algebra of operator fields  (see the definition of full algebra of
operator fields in 
[\textbf{5}],  [\textbf{8}] or [\textbf{17}]).  Also, one 
needs to use the fact that any element of a unital $C^*$-algebra that is 
sufficiently
close to the unit is invertible. 
Finally, to make $f_{i,j}$ uniformly within $\epsilon$ of $a_j$ (for all $i$,
$j$),
one needs the maximality of the algebra of operator fields $\A$ (see
[\textbf{8}] Proposition 1 and [\textbf{17}] Theorem 1.1).
     
    For simplicity, we will assume that there are only two intervals
($I_1$ and $I_2$) and only two $M$-tuples ( $(f_{1,1}, f_{1,2}, ..., f_{1,M})$
and $(f_{2,1}, f_{2,2}, ... , f_{2,M})$ ).  We may additionally assume
that neither interval is contained in the closure of the other, and that 
their intersection is a continuum.  Our goal is to ``connect" 
the two $M$-tuples over the intersection 
$\overline{ I_1 \bigcap I_2 }$ to get an element of
$Lg_{M} ( \A ) $ which approximates $(a_1, a_2, ... ,a_M )$ within 
$\epsilon$.  The proof for more than two intervals is an iteration of this
procedure (after appropriate contraction, removal or addition of 
intervals...). 

     Our procedure for ``connecting" the two $M$-tuples over $\overline{ I_1
\bigcap I_2}$ will involve constructing a sequence of strictly increasing
points $\{ t_n \}_{n=1}^{\infty}$ in $I_1 \bigcap I_2$, and constructing
sequences of operator fields
$\{ \alpha_{j}^n \}_{n = 1}^{\infty}$ ($j = 1, 2, ... , M$) satisfying:
\begin{enumerate} 
\item $\alpha_{j}^n$ is a continuous operator field in $\A(n)$, for all $j$,
$n$. 
Here 
$\A(n)$ is the unital maximal full algebra of operator fields gotten by 
restricting $\A$ to the interval $[t_n, t_{n+1}]$.  (In particular, the 
continuity structure for $\A(n)$ is gotten by taking the restriction to 
$[t_n, t_{n+1}]$ of all the fields in $\F$). 
\item  $\alpha_{j}^n$ is within $\epsilon$ of the restricion of $a_j$ to 
$[t_n, t_{n+1}]$, for all $j$, $n$.   
\item $\alpha_{j}^1 (t_1) = f_{1,j}(t_1)$ and $\alpha_{j}^{n+1} (t_{n+1} ) = 
\alpha_{j}^n (t_{n+1})$, for all $j$, $n$. 
\item $\alpha_{j}^{n+1}$ is within $\epsilon / {2^n}$ of the restriction 
of $f_{2,j}$ to $[t_{n+1}, t_{n+2}]$, for all $j$, $n$.  And 
\item $( \alpha_{1}^n, \alpha_{2}^n, ... , \alpha_{M}^n )$ is in 
$Lg_M ( \A(n) )$;  or equivalently, 
$\sum_{j = 1}^M \alpha_{j}^n (t)^* \alpha_{j}^n (t)
$ is an invertible element of $\A_t$ for all $t \in [t_n, t_{n+1}]$,  for
all $n$ (the proof of equivalence 
is a small spectral theory argument which uses the fact that 
$[t_n, t_{n+1}]$ is
compact).
\end{enumerate}
\noindent Henceforth, we will let ``(*)" denote conditions (1) - (5).

     By [\textbf{10}] Lemma 2.4 and our hypothesis for $M$, we have 
that $Lg_M ( \A_t ) \bigcap
\{ (b_1, b_2, ... , b_M ) \in \A_{t}^M : \| a_i (t) - b_i \| < \epsilon
\makebox{   } \forall i \}$  is a nonempty connected open set for all  
$t \in [0,1]$.  Hence, fixing $t_1 \in I_1 \bigcap I_2$, there are continuous
paths $\gamma_{ j} : [0,1] \rightarrow Lg_M ( \A_{t_1} )$, $j = 1,2, ...
, M$, such that a)  $\gamma_{j}(0) = f_{1,j}(t_1)$ and 
$\gamma_{j}(1) = f_{2, j}(t_1)$ for all $j$, and b) $\gamma_{j}(s)$
is within $\epsilon$ of $a_j (t_1)$ for all  $s \in [0,1]$ and for all $j$. 

     Note that since $[0,1]$ is compact, $( \gamma_{1}, 
\gamma_{2} ,...,
\gamma_{M} )$ is in $Lg_M ( C[0,1] \otimes \A_{t_1} )$.  Now consider
the unital maximal full algebra of operator fields $\widetilde{ \A }$, which 
has base space $[0,1]$ and fibre algebras $\{ C[0,1] \otimes \A_t \}_{ t \in 
[0,1] }$.  The continuity structure $\widetilde{ \F }$ for $\widetilde{ \A }$ 
consists of all operator fields of the form 
$ t \mapsto \sum_{i = 1}^N f_{i} \otimes c_{i} (t)$, $t \in [0,1]$.  Here
the $f_i$s are in $C[0,1]$, the $c_{i}$s are continuous operator fields in 
$\A$ (with respect to the continuity structure $\F$), and $N$ is a nonnegative
integer.  Now since $\widetilde{ \A }$ is a full algebra of operator fields,
there are continuous operator fields $t \mapsto \gamma_j (., t)$ in 
$\widetilde{\A}$ such that $\gamma_j (., t_1) = \gamma_j (.)$, $j= 1, 2, ... ,
M$.  (Here, the second variable of $\gamma_j (.,.)$ ranges over the base 
space of $\widetilde{\A}$, and for $t \in [0,1]$, $\gamma_j (.,t)$ is an    
element of the fibre algebra $C[0,1] \otimes \A_t$ of $\widetilde{\A}$).

     One can show that  
the operator field $t \mapsto \gamma_j (1,t)$ is a 
continuous operator field in $\A$ (with continuity structure $\F$) for 
$j = 1,.., M$.  Also, one may view the field $t \mapsto a_j(t)$ as an operator
field in $\widetilde{\A}$ in the natural way for $j=1,..,M$.  From these 
and the fact that (as elements of $C[0,1] \otimes \A_{t_1}$) 
$\gamma_j (., t_1)$ is within $\epsilon$ of $a_j (t_1)$, we can find an open
neighbourhood $V$ of $t_1$ with $V \subset I_1 \bigcap I_2$ such that 
for all $t \in V$, a)  as elements of $C[0,1] \otimes \A_t$,    
$\gamma_j (., t)$ is within $\epsilon$ of $a_j (t)$, $j=1,2,...,M$, and 
b)  $\gamma_j (1,t)$ is within $\epsilon/2$ of $f_{2,j} (t)$.   Furthermore, by 
continuity and since elements in a unital $C^*$-algebra which are sufficiently
close to 
the unit are invertible, we may assume that $V$ is sufficiently small so that
for all $t \in V$, $\sum_{j=1}^M \gamma_{j}(.,t)^* \gamma_{j}(., t)$ is an
invertible element of $C[0,1] \otimes \A_t$.

       One can show that  for any continuous function 
$g: [0,1] \rightarrow [0,1]$,
the operator field $t \mapsto \gamma_j ( g(t), t )$ is a continuous field in
$\A$, for $j = 1, 2, ..., M$.  Hence, we can pick $t_2 \in V$ such that 
$t_1 < t_2$;  and we can let $\alpha_{j}^1 (t) =_{df} \gamma_{j} (
(t - t_1 ) / ( t_2 - t_1 ) )$ for all $t \in [ t_1, t_2 ]$, $j= 1,2,..,M$.
This will give us the first members of the sequences in (*).

     Now we can repeat almost exactly the same argument as before, replacing
$t_1$ by the point $t_2$ and replacing $f_{1,j}(t_1)$ by $\alpha_{j}^1 (t_2)$,
$j = 1, 2, ..., M$.  In this manner, we get a point $t_3$ and operator fields
$\alpha_{j}^2$, $j = 1, 2, ..., M$, which will be the next members of the 
sequences in (*).  We need only note two minor modifications that are needed in
the argument:  
a)  firstly, one has to use the fact that for all $t$, the 
set \[
Lg_M ( \A_{t} ) \bigcap \{ (b_1, b_2, ..., b_M ) \in {\A_{ t}}^M :
\| a_j (t ) - b_j \| < \epsilon \makebox{   and   } 
\| f_{2,j} (t) - b_j \| < \epsilon/2,  \makebox{    for all   } j \}
\]
\noindent is a connected open set (see [\textbf{10}] Lemma 2.4), 
and b) when choosing the        
corresponding neighbourhood about $t_2$, one must make it sufficiently small
so that the corresponding quantities which result will also be sufficiently
small in order to fulfill condition (4) in (*). 

     Repeating this process ad infinitum (making the appropriate
modifications at each step), we get a sequence of points 
$\{ t_n \}_{n=1}^{\infty}$ and sequences of continuous operator fields
$\{ \alpha_{j}^n \}_{ n = 1}^{\infty}$, $j=1, 2, ..., M$, which fulfill the 
conditions in (*).  Now let $\tilde{t} = lim_{n \rightarrow \infty} t_n$.
For $j= 1, 2,..., M$, let $\alpha_j$ be the continuous operator field in
$\A$ defined by 
\begin{enumerate}
\item $\alpha_{j}(t) = f_{1,j} (t)$ for $t \in [0, t_1]$,
\item $\alpha_{j}(t) = \alpha_{j}^{n} (t)$ for $t \in [t_n, t_{n+1} ]$, and
\item $\alpha_{j} (t) = f_{2,j}(t)$ for all $t \in [ \tilde{t}, 1]$.
\end{enumerate}  
\noindent Then $( \alpha_1, \alpha_2, ..., \alpha_M) \in Lg_M ( \A )$, and 
for all $j = 1,2, ..., M$, $\alpha_j$ approximates $a_j$ within $\epsilon$.
\end{proof}

\begin{proof}[Proof of Theorem 1.1] 
We proceed by induction.  The base case $k=1$ has already been dealt with 
in Lemma 2.1.  We now do the induction step, supposing that $k \geq 2$. 
By [\textbf{8}] Theorem 4, let $\pi : Prim (\A) \rightarrow [0,1]^k$ 
be the continuous open 
surjection corresponding to the continuous field decomposition of $\A$ in 
the hypothesis (Here $Prim(\A)$ is the primitive ideal space of $\A$).   
 
    Since $k \geq 2$, the map $p : [0,1]^k          
\rightarrow [0,1]^{k-1}$, given by projecting onto the first $k-1$ coordinates,
is a continuous open surjection.  Hence, the composition 
$p \circ \pi : Prim( \A ) \rightarrow [0,1]^{k-1}$ is a continuous open 
surjection.  Hence by [\textbf{8}] Theorem 4, we can realize $\A$ as a 
unital maximal full 
algebra of operator fields with base space $[0,1]^{k-1}$  and fibre algebras,
say, $\{ \B_s \}_{ s \in [0,1]^{k-1} }$.  
Hence by the induction hypothesis,
$sr ( \A ) \leq sup_{ s \in [0,1] } sr ( C([0,1]^{k-1}) \otimes \B_s )$.  
                       
    But for each $s$, the fibre algebra $\B_s$ can be realized 
as a maximal full algebra
of operator fields with base space $p^{-1}(s) = \{ s \} \times [0,1]$ and 
fibre algebras $\{ \A_r \}_{ r \in \{ s \} \times [0,1] }$ (By [\textbf{8}] Theorem 
4,  
$\B_s$ is isomorphic to  
$\A / I$ where 
$I = \bigcap (p \circ \pi)^{-1} (s)$.  Using this fact, one can construct the
natural continuous open surjection of $Prim( \B_ s)$ 
onto $ \{ s \} \times [0,1]$).
Let us suppose that this continuous field decomposition of $\B_s$ is given
by a  continuity structure $\mathcal{G}$.  Then $C([0,1]^{k-1}) \otimes \B_s$ 
can 
be realized as a unital maximal full algebra of operator fields with
base space $[0,1]$ and fibre algebras $\{ C([0,1]^{k-1}) \otimes \A_r \}_{
r \in \{ s \} \times [0,1] }$.  Here the continuity structure consists of
operator fields of the form $ r \mapsto 
\sum_{i = 1}^N f_{i} \otimes b_{i} (r)$, for $r \in \{ s \} \times [0,1]$.
Here the $f_i$s are in $C([0,1]^{k - 1})$, the $b_{i}$s are continuous fields
in $\B_s$ (with respect to the continuity structure $\mathcal{G}$) and 
$N$ is a     
nonnegative integer.  

     Hence by the induction hypothesis, 
for each $s$, we have that $sr ( C([0,1]^{k-1}) \otimes 
\B_s ) \leq sup_{r \in \{ s \} \times [0,1] } sr( C([0,1]^k ) \otimes \A_r )$.
It follows, then, that $sr( \A ) \leq sup_{ t \in [0,1]^k } sr ( C([0,1]^{k})
\otimes \A_t )$.
\end{proof}

     We note that the statements of Lemma 2.1  and Theorem 1.1
 would still be true if the
unit interval $[0,1]$ was replaced by the circle $S^1$ (and if the $k$-cube
$[0,1]^k$ was replaced by the $k$-torus $\mathbb{T}^k$).  The proofs would
be exactly the same.  In other words, we have that
  
\begin{cor}  Let $\A$ be a unital maximal full algebra of operator fields
with base space the $k$-torus $\T^k$ and fibre algebras $\{ \A_t \}_{t \in
\T^k}$.  Then the stable rank of $\A$ satisfies
\[
sr( \A ) \leq sup_{ t \in \T^k } sr( C([0,1]^k ) \otimes \A_t ). 
\]
\end{cor}
 
\begin{thm}  Let $H_{2n +1}^{\Z}$ be the discrete Heisenberg group of 
rank $2n + 1$.  Let $C^*( H_{2n +1}^{\Z} )$ be the univeral $C^*$-algebra of
$H_{2n +1}^{\Z}$.  Then $sr( C^*( H_{2n +1}^{\Z} ) ) = n+1$.
\end{thm} 
\begin{proof}
By [\textbf{12}] Proposition 1.7 and Theorem 4.3, 
 and the fact that $C( \T^{2n} )$ is a quotient of 
$C^*( H_{2n + 1}^{\Z} )$, the stable rank of $C^*( H_{2n + 1}^{\Z} )$ is  
greater than or equal to $n+1$.  
By [\textbf{1}] and [\textbf{9}] Theorem 3.4, 
$C^*( H_{2n + 1}^{\Z} )$ can be realized 
as a unital maximal full algebra of operator fields with base space 
the $1$-torus $\T$ and fibre algebras, say, $\{ \A_t \}_{ t \in \T}$.         
Hence, by the Corollary 2.2, 
the stable rank of $C^*( H_{2n + 1}^{\Z} )$ satisfies 
$sr ( C^*( H_{2n + 1}^{\Z} )) \leq sup_{ t \in \T } sr( C[0,1] \otimes \A_t )$.

     Now by [\textbf{1}] and 
[\textbf{9}] Theorem 3.4, each fibre algebra 
$\A_t$ can be     
realized as a unital maximal full algebra of operator fields with base space
a torus with dimension less than or equal to $2n$ (the zero-dimensional torus
being a point) and fibre algebras, 
say, $\{ B_{s}^t \}_{ s \in \T^{l_t} }$ where
$\T^{l_t}$ is the base in the continuous field decomposition of $\A_t$ ($l_t 
\leq
2n$).  Moreover, for each $t \in \T$ and for each
$s \in \T^{l_t}$, $\B_{s}^t$ is isomorphic to either
a full matrix algebra $\M_n (\C)$ or $\M_n ( \C ) \otimes (\bigotimes^n  
\mathbb{A}_{\theta})$, where $\bigotimes^n \mathbb{A}_{\theta}$ is the $n$ 
times tensor product
of a fixed irrational rotation algebra $\mathbb{A}_{\theta}$ with 
irrational rotation
$\theta$.  Now by [\textbf{4}] each irrational rotation algebra can be decomposed
as an inductive limit of building blocks of the form $M_m (C( \T )) \oplus
\M_n (C (\T))$ (the integers $m$ and $n$ get arbitrarily large as we move
along building blocks in the inductive limit).  Hence, since $C[0,1] \otimes
\A_t$ can be realized as a maximal full algebra of operator fields with
base space $\T^{l_t}$ and fibre algebras 
$\{ C[0,1] \otimes \B_{s}^t \}_{ s \in \T}$,
it follows, by Corollary 2.2, [\textbf{12}] Proposition 1.7 and 
Theorems 5.1 and 
6.1, 
that 
$sr( C[0,1] \otimes \A_t ) \leq n+1$, for all $t$.  Hence, 
$sr( C^*( H_{2n + 1}^{\Z} ) = n+1$.     
\end{proof}

\section*{References} 

\newcounter{bean}
\begin{list}{\textbf{\arabic{bean}}}{\usecounter{bean}}             
\item   \textsc{J. Anderson and W. Paschke}, 
`The rotation algebra',
\emph{Houston J. Math.}, 15, (1989), 1 - 26. 
\item \textsc{L. Baggett and J. Packer}, 
`The primitive ideal space of two-step nilpotent group $C^*$-algebras',
\emph{J. Funct. Anal.}, 124, (1994), 389-426.

\item \textsc{ K. R. Davidson},
\emph{$C^*$-algebras by example},
Fields Institute Monographs, 6.

\item \textsc{G. A. Elliott and D. E. Evans}, `The structure of
the irrational rotation $C^*$-algebra', \emph{Ann. of Math. (2)}, 138,
no. 3, (1993), 477-501.

\item \textsc{J. M. G. Fell},
`The structure of algebras of operator fields', \emph{Acta Math.}, 
106, (1961), 233-280.

\item \textsc{N. E. Hassan}, `Rangs Stables de certaines 
extensions', \emph{J. London Math. Soc.}, (2), 52, (1995), 605-624.

\item \textsc{R. Herman, L. N. Vaserstein}, 
`The stable rank of $C^*$-algebras', \emph{Invent. Math.}, 
77, (1984), 553-555.

\item \textsc{R-Y Lee}, `On $C^*$-algebras of operator fields',
\emph{Indiana U. Math. Journal}, 25, no. 4, (1976), 303-314.  

\item  \textsc{S. T. Lee and J. Packer},
`Twisted group $C^*$-algebras for two-step nilpotent and generalized
discrete Heisenberg group', \emph{J. Operator Th.}, 34, no. 1, (1995), 
91-124.

\item \textsc{V. Nistor},
`Stable range for the tensor products of extensions of $K$ by 
$C(X)$', \emph{J. Operator Th.}, 16, (1986), 387-396.

\item \textsc{J. Packer and I. Raeburn}, 
`On the structure of twisted group $C^*$-algebras',
\emph{Trans. AMS}, 334, no. 2, (1992), 685-718.

\item \textsc{M. A. Rieffel},
`Dimension and stable rank in the $K$-theory of $C^*$-algebras',
\emph{Proc. London Math. Soc.}, (3), 46, (1983), 301-333.

\item \textsc{T. Sudo}, `Dimension theory of group $C^*$-algebras of 
connected Lie groups of type I',  \emph{J. Math. Soc. Japan}, 52, no. 3,
(2000), 583-590.    

\item \textsc{T. Sudo}, 'Structure of group $C^*$-algebras of Lie semi-direct
products $\mathbb{C} \times | \mathbb{R}$', 
\emph{J. Operator Theory}, 46, (2001), 25-38.

\item \textsc{T. Sudo and H. Takai}, `Stable rank of the 
$C^*$-algebras of nilpotent Lie groups', \emph{Internat. J. Math.}, 6,
no. 3, (1995), 439-446.

\item \textsc{T. Sudo and H. Takai}, `Stable rank of the 
$C^*$-algebras of solvable Lie groups of type I', \emph{J. Operator Theory},
38, no. 1, (1997), 67-86.

\item \textsc{J. Tomiyama}, 
`Topological representation of $C^*$-algebras',
\emph{Tohoku Math. J.}, (2), 14, (1962), 187-204.
\end{list}

\end{document}